\expandafter\chardef\csname pre amssym.def at\endcsname=\the\catcode`\@ 
\catcode`\@=11 
 
\def\undefine#1{\let#1\undefined} 
\def\newsymbol#1#2#3#4#5{\let\next@\relax 
 \ifnum#2=\@ne\let\next@\msafam@\else 
 \ifnum#2=\tw@\let\next@\msbfam@\fi\fi 
 \mathchardef#1="#3\next@#4#5} 
\def\mathhexbox@#1#2#3{\relax 
 \ifmmode\mathpalette{}{\m@th\mathchar"#1#2#3}%
 \else\leavevmode\hbox{$\m@th\mathchar"#1#2#3$}\fi} 
\def\hexnumber@#1{\ifcase#1 0\or 1\or 2\or 3\or 4\or 5\or 6\or 7\or 8\or 
 9\or A\or B\or C\or D\or E\or F\fi} 
 
\font\tenmsa=msam10 
\font\sevenmsa=msam7 
\font\fivemsa=msam5 
\newfam\msafam 
\textfont\msafam=\tenmsa 
\scriptfont\msafam=\sevenmsa 
\scriptscriptfont\msafam=\fivemsa 
\edef\msafam@{\hexnumber@\msafam} 
\mathchardef\dabar@"0\msafam@39 
\def\dashrightarrow{\mathrel{\dabar@\dabar@\mathchar"0\msafam@4B}} 
\def\dashleftarrow{\mathrel{\mathchar"0\msafam@4C\dabar@\dabar@}} 
 
\def\ulcorner{\delimiter"4\msafam@70\msafam@70 } 
\def\urcorner{\delimiter"5\msafam@71\msafam@71 } 
\def\llcorner{\delimiter"4\msafam@78\msafam@78 } 
\def\lrcorner{\delimiter"5\msafam@79\msafam@79 } 
\def\yen{{\mathhexbox@\msafam@55 }} 
\def\checkmark{{\mathhexbox@\msafam@58 }} 
\def\circledR{{\mathhexbox@\msafam@72 }} 
\def\maltese{{\mathhexbox@\msafam@7A }} 
 
\font\tenmsb=msbm10 
\font\sevenmsb=msbm7 
\font\fivemsb=msbm5 
\newfam\msbfam 
\textfont\msbfam=\tenmsb 
\scriptfont\msbfam=\sevenmsb 
\scriptscriptfont\msbfam=\fivemsb 
\edef\msbfam@{\hexnumber@\msbfam}

\catcode`\@=\csname pre amssym.def at\endcsname 
 
\expandafter\ifx\csname pre amssym.tex at\endcsname\relax \else \endinput\fi 
\expandafter\chardef\csname pre amssym.tex at\endcsname=\the\catcode`\@ 
\catcode`\@=11 
\newsymbol\boxdot 1200 
\newsymbol\boxplus 1201 
\newsymbol\boxtimes 1202 
\newsymbol\square 1003 
\newsymbol\blacksquare 1004 
\newsymbol\centerdot 1205 
\newsymbol\lozenge 1006 
\newsymbol\blacklozenge 1007 
\newsymbol\circlearrowright 1308 
\newsymbol\circlearrowleft 1309 
\undefine\rightleftharpoons 
\newsymbol\rightleftharpoons 130A 
\newsymbol\leftrightharpoons 130B 
\newsymbol\boxminus 120C 
\newsymbol\Vdash 130D 
\newsymbol\Vvdash 130E 
\newsymbol\vDash 130F 
\newsymbol\twoheadrightarrow 1310 
\newsymbol\twoheadleftarrow 1311 
\newsymbol\leftleftarrows 1312 
\newsymbol\rightrightarrows 1313 
\newsymbol\upuparrows 1314 
\newsymbol\downdownarrows 1315 
\newsymbol\upharpoonright 1316 
  
\newsymbol\downharpoonright 1317 
\newsymbol\upharpoonleft 1318 
\newsymbol\downharpoonleft 1319 
\newsymbol\rightarrowtail 131A 
\newsymbol\leftarrowtail 131B 
\newsymbol\leftrightarrows 131C 
\newsymbol\rightleftarrows 131D 
\newsymbol\Lsh 131E 
\newsymbol\Rsh 131F 
\newsymbol\rightsquigarrow 1320 
\newsymbol\leftrightsquigarrow 1321 
\newsymbol\looparrowleft 1322 
\newsymbol\looparrowright 1323 
\newsymbol\circeq 1324 
\newsymbol\succsim 1325 
\newsymbol\gtrsim 1326 
\newsymbol\gtrapprox 1327 
\newsymbol\multimap 1328 
\newsymbol\therefore 1329 
\newsymbol\because 132A 
\newsymbol\doteqdot 132B 
  
\newsymbol\triangleq 132C 
\newsymbol\precsim 132D 
\newsymbol\lesssim 132E 
\newsymbol\lessapprox 132F 
\newsymbol\eqslantless 1330 
\newsymbol\eqslantgtr 1331 
\newsymbol\curlyeqprec 1332 
\newsymbol\curlyeqsucc 1333 
\newsymbol\preccurlyeq 1334 
\newsymbol\leqq 1335 
\newsymbol\leqslant 1336 
\newsymbol\lessgtr 1337 
\newsymbol\backprime 1038 
\newsymbol\risingdotseq 133A 
\newsymbol\fallingdotseq 133B 
\newsymbol\succcurlyeq 133C 
\newsymbol\geqq 133D 
\newsymbol\geqslant 133E 
\newsymbol\gtrless 133F 
\newsymbol\sqsubset 1340 
\newsymbol\sqsupset 1341 
\newsymbol\vartriangleright 1342 
\newsymbol\vartriangleleft 1343 
\newsymbol\trianglerighteq 1344 
\newsymbol\trianglelefteq 1345 
\newsymbol\bigstar 1046 
\newsymbol\between 1347 
\newsymbol\blacktriangledown 1048 
\newsymbol\blacktriangleright 1349 
\newsymbol\blacktriangleleft 134A 
\newsymbol\vartriangle 134D 
\newsymbol\blacktriangle 104E 
\newsymbol\triangledown 104F 
\newsymbol\eqcirc 1350 
\newsymbol\lesseqgtr 1351 
\newsymbol\gtreqless 1352 
\newsymbol\lesseqqgtr 1353 
\newsymbol\gtreqqless 1354 
\newsymbol\Rrightarrow 1356 
\newsymbol\Lleftarrow 1357 
\newsymbol\veebar 1259 
\newsymbol\barwedge 125A 
\newsymbol\doublebarwedge 125B 
\undefine\angle 
\newsymbol\angle 105C 
\newsymbol\measuredangle 105D 
\newsymbol\sphericalangle 105E 
\newsymbol\varpropto 135F 
\newsymbol\smallsmile 1360 
\newsymbol\smallfrown 1361 
\newsymbol\Subset 1362 
\newsymbol\Supset 1363 
\newsymbol\Cup 1264 
  
\newsymbol\Cap 1265 
  
\newsymbol\curlywedge 1266 
\newsymbol\curlyvee 1267 
\newsymbol\leftthreetimes 1268 
\newsymbol\rightthreetimes 1269 
\newsymbol\subseteqq 136A 
\newsymbol\supseteqq 136B 
\newsymbol\bumpeq 136C 
\newsymbol\Bumpeq 136D 
\newsymbol\lll 136E 
  
\newsymbol\ggg 136F 
  
\newsymbol\circledS 1073 
\newsymbol\pitchfork 1374 
\newsymbol\dotplus 1275 
\newsymbol\backsim 1376 
\newsymbol\backsimeq 1377 
\newsymbol\complement 107B 
\newsymbol\intercal 127C 
\newsymbol\circledcirc 127D 
\newsymbol\circledast 127E 
\newsymbol\circleddash 127F 
\newsymbol\lvertneqq 2300 
\newsymbol\gvertneqq 2301 
\newsymbol\nleq 2302 
\newsymbol\ngeq 2303 
\newsymbol\nless 2304 
\newsymbol\ngtr 2305 
\newsymbol\nprec 2306 
\newsymbol\nsucc 2307 
\newsymbol\lneqq 2308 
\newsymbol\gneqq 2309 
\newsymbol\nleqslant 230A 
\newsymbol\ngeqslant 230B 
\newsymbol\lneq 230C 
\newsymbol\gneq 230D 
\newsymbol\npreceq 230E 
\newsymbol\nsucceq 230F 
\newsymbol\precnsim 2310 
\newsymbol\succnsim 2311 
\newsymbol\lnsim 2312 
\newsymbol\gnsim 2313 
\newsymbol\nleqq 2314 
\newsymbol\ngeqq 2315 
\newsymbol\precneqq 2316 
\newsymbol\succneqq 2317 
\newsymbol\precnapprox 2318 
\newsymbol\succnapprox 2319 
\newsymbol\lnapprox 231A 
\newsymbol\gnapprox 231B 
\newsymbol\nsim 231C 
\newsymbol\ncong 231D 
\newsymbol\diagup 231E 
\newsymbol\diagdown 231F 
\newsymbol\varsubsetneq 2320 
\newsymbol\varsupsetneq 2321 
\newsymbol\nsubseteqq 2322 
\newsymbol\nsupseteqq 2323 
\newsymbol\subsetneqq 2324 
\newsymbol\supsetneqq 2325 
\newsymbol\varsubsetneqq 2326 
\newsymbol\varsupsetneqq 2327 
\newsymbol\subsetneq 2328 
\newsymbol\supsetneq 2329 
\newsymbol\nsubseteq 232A 
\newsymbol\nsupseteq 232B 
\newsymbol\nparallel 232C 
\newsymbol\nmid 232D 
\newsymbol\nshortmid 232E 
\newsymbol\nshortparallel 232F 
\newsymbol\nvdash 2330 
\newsymbol\nVdash 2331 
\newsymbol\nvDash 2332 
\newsymbol\nVDash 2333 
\newsymbol\ntrianglerighteq 2334 
\newsymbol\ntrianglelefteq 2335 
\newsymbol\ntriangleleft 2336 
\newsymbol\ntriangleright 2337 
\newsymbol\nleftarrow 2338 
\newsymbol\nrightarrow 2339 
\newsymbol\nLeftarrow 233A 
\newsymbol\nRightarrow 233B 
\newsymbol\nLeftrightarrow 233C 
\newsymbol\nleftrightarrow 233D 
\newsymbol\divideontimes 223E 
\newsymbol\varnothing 203F 
\newsymbol\nexists 2040 
\newsymbol\Finv 2060 
\newsymbol\Game 2061 
\newsymbol\mho 2066 
\newsymbol\eth 2067 
\newsymbol\eqsim 2368 
\newsymbol\beth 2069 
\newsymbol\gimel 206A 
\newsymbol\daleth 206B 
\newsymbol\lessdot 236C 
\newsymbol\gtrdot 236D 
\newsymbol\ltimes 226E 
\newsymbol\rtimes 226F 
\newsymbol\shortmid 2370 
\newsymbol\shortparallel 2371 
\newsymbol\smallsetminus 2272 
\newsymbol\thicksim 2373 
\newsymbol\thickapprox 2374 
\newsymbol\approxeq 2375 
\newsymbol\succapprox 2376 
\newsymbol\precapprox 2377 
\newsymbol\curvearrowleft 2378 
\newsymbol\curvearrowright 2379 
\newsymbol\digamma 207A 
\newsymbol\varkappa 207B 
\newsymbol\Bbbk 207C 
\newsymbol\hslash 207D 
\undefine\hbar 
\newsymbol\hbar 207E 
\newsymbol\backepsilon 237F 
\catcode`\@=\csname pre amssym.tex at\endcsname 
 
\magnification=1200 
\hsize=468truept 
\vsize=646truept 
\voffset=-10pt 
\parskip=1pc 
\baselineskip=14truept 
\count0=1 
 
\dimen100=\hsize 
 
\def\leftill#1#2#3#4{ 
\medskip 
\line{$ 
\vcenter{ 
\hsize = #1truept \hrule\hbox{\vrule\hbox to  \hsize{\hss \vbox{\vskip#2truept 
\hbox{{\copy100 \the\count105}: #3}\vskip2truept}\hss } 
\vrule}\hrule} 
\dimen110=\dimen100 
\advance\dimen110 by -36truept 
\advance\dimen110 by -#1truept 
\hss \vcenter{\hsize = \dimen110 
\medskip 
\noindent { #4\par\medskip}}$} 
\advance\count105 by 1 
} 
\def\rightill#1#2#3#4{ 
\medskip 
\line{ 
\dimen110=\dimen100 
\advance\dimen110 by -36truept 
\advance\dimen110 by -#1truept 
$\vcenter{\hsize = \dimen110 
\medskip 
\noindent { #4\par\medskip}} 
\hss \vcenter{ 
\hsize = #1truept \hrule\hbox{\vrule\hbox to  \hsize{\hss \vbox{\vskip#2truept 
\hbox{{\copy100 \the\count105}: #3}\vskip2truept}\hss } 
\vrule}\hrule} 
$} 
\advance\count105 by 1 
} 
\def\midill#1#2#3{\medskip 
\line{$\hss 
\vcenter{ 
\hsize = #1truept \hrule\hbox{\vrule\hbox to  \hsize{\hss \vbox{\vskip#2truept 
\hbox{{\copy100 \the\count105}: #3}\vskip2truept}\hss } 
\vrule}\hrule} 
\dimen110=\dimen100 
\advance\dimen110 by -36truept 
\advance\dimen110 by -#1truept 
\hss $} 
\advance\count105 by 1 
} 
\def\insectnum{\copy110\the\count120 
\advance\count120 by 1 
}

\font\ninerm=cmr9 
\font\eightrm=cmr8

\font\tenrm=cmr10 at 10pt 
 
\font\sc=cmcsc10

 
\def\msb{\fam\msbfam\tenmsb} 
 
\def\bba{{\msb A}} 
 
\def\bbc{{\msb C}}

\def\bbp{{\msb P}} 
\def\bbq{{\msb Q}} 
\def\bbr{{\msb R}}

\def\bbz{{\msb Z}}

\def\grD{\Delta}

\def\grL{\Lambda} 
\def\grO{\Omega}

\def\gra{\alpha}

\def\grd{\delta} 
\def\gre{\epsilon}

\def\gri{\iota}

\def\gro{\omega}

\def\grr{\rho}

\def\grz{\zeta} 
 
\def\la#1{\hbox to #1pc{\leftarrowfill}} 
\def\ra#1{\hbox to #1pc{\rightarrowfill}} 
 
\def\fract#1#2{\raise4pt\hbox{$ #1 \atop #2 $}} 
\def\decdnar#1{\phantom{\hbox{$\scriptstyle{#1}$}} 
\left\downarrow\vbox{\vskip15pt\hbox{$\scriptstyle{#1}$}}\right.} 
 
\def\bowtie{\hbox to 1pt{\hss}\raise.66pt\hbox{$\scriptstyle{>}$} 
\kern-4.9pt\triangleleft} 
\def\hsmash{\triangleright\kern-4.4pt\raise.66pt\hbox{$\scriptstyle{<}$}} 
\def\boxit#1{\vbox{\hrule\hbox{\vrule\kern3pt 
\vbox{\kern3pt#1\kern3pt}\kern3pt\vrule}\hrule}}

\def\za{\vrule height6pt width4pt depth1pt}

\font\aa=eufm10

\def\Got#1{\hbox{\aa#1}}

\def\bfw{{\bf w}}

\def\calo{{\cal O}} 
 
\def\cald{{\cal D}} 
 
\def\calf{{\cal F}}

\def\calo{{\cal O}}

\def\cals{{\cal S}}

\def\calz{{\cal Z}}

\def\gF{{\Got F}} 
\def\gG{{\Got G}}

\def\Got#1{\hbox{\aa#1}}

\def\gsp1{{\Got s}{\Got p}(1)}

\font\svtnrm=cmr17

\font\bsc=cmcsc10 at 10truept

\def\Ric{\hbox{Ric}}

\def\Se{Sasakian-Einstein }

\def\inv{1}
\def\pos{2}
\def\consum{3}

\phantom{ooo}
\bigskip\bigskip
\centerline{\svtnrm On Positive Sasakian Geometry}  \medskip

\bigskip\bigskip
\centerline{\sc Charles P. Boyer~~ Krzysztof Galicki~~ Michael Nakamaye}
\footnote{}{\ninerm During the preparation of this work the first two authors 
were partially supported by NSF grant DMS-9970904, and third author by NSF
grant DMS-0070190.} \bigskip

\centerline{\vbox{\hsize = 5.85truein
\baselineskip = 12.5truept
\eightrm
\noindent {\bsc Abstract:}
A Sasakian structure $\scriptstyle{\cals=(\xi,\eta,\Phi,g)}$ on a manifold 
$\scriptstyle{M}$ is called {\it positive} if its basic first Chern
class $\scriptstyle{c_1(\calf_\xi)}$ can be represented by a positive
$\scriptstyle{(1,1)}$-form with respect to its transverse holomorphic
CR-structure. We prove a theorem that says that every positive Sasakian
structure can be deformed to a Sasakian structure whose metric has positive
Ricci curvature. This allows us by example to give a completely independent
proof of a result of Sha and Yang [SY] that for every positive integer
$\scriptstyle{k}$  the 5-manifolds $\scriptstyle{k\#(S^2\times S^3)}$
admit metrics of positive Ricci curvature.}} \tenrm

\vskip .5in
\baselineskip = 10 truept
\centerline{\bf Introduction}  
\bigskip

Let $(M,J)$ is a compact complex
manifold and $g$ a K\"ahler metric on $M$, with K\"ahler form M.
Suppose that $\rho'$ is a real, closed $(1,1)$-form on
$M$ with $[\rho']=2\pi c_1(M)$. Then there exists a unique K\"ahler
metric $g'$ on $M$ with K\"ahler form $\omega'$, such that
$[\omega]=[\omega']\in H^2(M,\bbr)$, and the Ricci form of $g'$ is
$\rho'$. The above statement is the celebrated Calabi Conjecture which
was posed by Eugene Calabi in 1954. The conjecture in its full generality was
eventually proved by Yau in 1976. In the Fano case when $c_1(M)>0$, i.e.,
when the first Chern class can be represented by a positive-definite 
real, closed $(1,1)$-form $\rho'$ on $M$, the conjecture implies that
the K\"ahler form of $M$ can be represented by a metric of
positive Ricci curvature. 

There are several reasons one might be interested in a more
general Calabi Problem when $M$ is not necessarily
a smooth manifold but rather a $V$-manifold or an orbifold [DK]. In the
context of Sasakian manifolds one would naturally consider the 
K\"ahler geometry of the associated one-dimensional foliation. In this
context one can actually prove a ``transverse Yau theorem" and this was done
by El Kacimi-Alaoui in 1990 [ElK]. In this note we discuss some
consequences of the transverse Yau theorem. In particular, we prove the
following which can be viewed as a Sasakian version of the ``positive Calabi
Conjecture'' mentioned above:

\medskip
\noindent{\sc Theorem A}: \tensl Let $\cals=(\xi,\eta,\Phi,g)$ be a positive
Sasakian structure on a compact manifold $M$ of dimension $2n+1.$ Then $M$
admits a Sasakian structure $\cals'=(\xi',\eta',\Phi',g')$  with positive
Ricci curvature homologous to $\cals.$
\tenrm
\medskip
Using some of the techniques developed earlier in our studies of
circle V-bundles over log del Pezzo surfaces [BG2,BGN1,BGN2] we are able to
apply Theorem A to demonstrate 
\medskip
\noindent{\sc Theorem B}: \tensl There exist Sasakian metrics
with positive Ricci tensor on $k\#(S^2\times S^3)$ for every positive integer
$k.$ In particular, for each $k>5$ there is a $2(k-5)$ parameter family of
inequivalent Sasakian structures with positive Ricci curvature.  \tenrm  

From the standpoint of Riemannian metrics with positive Ricci curvature, this
result is not new as it was first obtained by Sha and Yang [SY] in 1991;
however, the Sasakian nature of such metrics is new.  Moreover, our method of
proof here is completely different.

In [JK] a list of all anticanonically (i.e. orbifold Fano index $=1$) 
quasi-smooth log del Pezzo surfaces embedded in a weighted projective space
$\bbp(\bfw)$ was given, and in [BGN1] it is shown that the largest
Picard number occuring is 10. However, as soon as we allow higher Fano index,
namely 2, we find: 

\medskip
\noindent{\sc Corollary C}: \tensl There exist log del Pezzo surfaces $\subset
\bbp(\bfw)$ with arbitrary Picard number. \tenrm

\bigskip
\baselineskip = 10 truept
\centerline{\bf \inv. Some Transverse Holomorphic Invariants}
\bigskip

Recall that associated with every Sasakian structure $\cals=(\xi,\eta,\Phi,g)$
on a smooth manifold $M$ is  its characteristic foliation $\calf_\xi,$ and as 
described elsewhere (cf. [BGN1]) the transverse geometry of a Sasakian manifold
is K\"ahler. Moreover,  there are transverse versions of both the Hodge theory
and the Dolbeault theory [ElK] of a Sasakian manifold. In particular if
$(\cals,\xi,\eta,\Phi,g)$ is a compact Sasakian manifold, then $d\eta$ defines
a nontrivial class $[d\eta]_B$ in the basic cohomology group
$H^{1,1}_B(\calf_\xi)\subset H^2_B(\calf_\xi).$ We also define the set
$\gF(\xi)$ of all deformed Sasakian structures $(\xi,\eta',\Phi',g')$ on
$\cals$ that are {\it homologous} to $(\xi,\eta,\Phi,g),$ that is such that
$[d\eta']_B=[d\eta]_B.$ Now we have the {\it basic Betti numbers}  
$$b^B_r(\calf_\xi)=\hbox{dim}~H^r_B(\calf_\xi), \leqno{\inv.1}$$ 
and the {\it basic hodge numbers}
$$h^{p,q}_B(\calf_\xi)= \hbox{dim}~H^{p,q}_B(\calf_\xi) \leqno{\inv.2}$$
which satisfy the relations
$$b^B_r(\calf_\xi)=\sum_{p+q=r}h^{p,q}(\calf_\xi), \qquad
h^{p,q}_B(\calf_\xi)=h^{q,p}_B(\calf_\xi).\leqno{\inv.3}$$
Accordingly, we have the {\it transverse Euler characteristic} of
$\calf_\xi$ given by 
$$\chi(\calf_\xi)=\sum_{p=0}^{2n}
(-1)^{p}\hbox{dim}~H^p_B(\calf),\leqno{\inv.4a}$$
and the {\it transverse holomorphic Euler characteristic} defined by
$$\chi_{hol}(\calf_\xi)=\sum_i (-1)^i\hbox{dim}~H^{0,q}_B(\calf_\xi).
\leqno{\inv.4b}$$  

Now the basic cohomology ring $H^*(\calf_\xi)$ is invariant under foliated
homeomorphisms of $M,$ that is homeomorphisms that preserve the foliation
[ElKN], so the basic Betti numbers are also invariant under such
homeomorphisms. The basic Hodge numbers, however, are only invariant under
homemorphisms that preserve the foliation together with its transverse complex
structure.  Thus, they are not only invariants of the Sasakian structure, but
also of its basic deformation class $\gF(\xi).$  In the case that
$\cals=(\xi,\eta,\Phi,g)$ is quasi-regular, the basic Betti numbers and basic
Hodge numbers coincide with the corresponding Betti numbers and Hodge numbers
on the space of leaves $\calz= \cals/S^1$ which is a compact K\"ahler
orbifold. 

The relationship between the basic cohomology $H^r_B(\calf_\xi)$ and the
ordinary cohomology $H^r(\cals,\bbr)$ is given by the exact sequence [Ton]
$$\cdots\ra{2.5}H_B^p(\calf_\xi) \fract{\gri_*}{\ra{2.5}}H^p(\cals,\bbr)
\fract{j_p}{\ra{2.5}} H_B^{p-1}(\calf_\xi) \fract{\grd}{\ra{2.5}}
H^{p+1}_B(\calf_\xi)\ra{2.5}\cdots   \leqno{\inv.5}$$  
where $\grd$ is the connecting homomorphism given by $\grd[\gra]_B=[d\eta\wedge
\gra]_B=[d\eta]_B\cup[\gra]_B,$ and $j_p$ is the composition of the map induced
by $\xi\rfloor$ with the well known isomorphism $H^r(M,\bbr)\approx
H^r(M,\bbr)^{S^1}$ where $H^r(M,\bbr)^{S^1}$ is the $S^1$-invariant
cohomology defined from the  $S^1$-invariant r-forms $\grO^r(M)^{S^1}.$
Another important transverse invariant of the class $\gF(\xi)$ is the basic
first Chern class $c_1(\calf_\xi)\in H^2(\calf_\xi)$ of the foliation
$\calf_\xi.$ The image $\gri_*(c_1(\calf_\xi))$ is the real first Chern class
of the normal bundle $\nu(\calf_\xi)$ to the foliation with its induced
transverse complex structure.

If the dimension of $\cals$ is $4n+1$ we can also define the {\it basic
Hirzebruch signature} $\tau_B(\calf_\xi)$ to be the signature of the bilinear
form defined by the cup product on the middle basic cohomology group
$H^{2n}_B(\calf_\xi).$ Using transverse Lefschetz theory [ElK] one can obtain
the usual formula:
$$\tau_B(\calf_\xi)=\sum_{p,q}(-1)^ph^{p,q}_B(\calf_\xi)=\sum_{p\equiv
q(2)}(-1)^ph^{p,q}_B(\calf_\xi). \leqno{\inv.6}$$

As is usual in complex geometry we introduce the {\it geometric genus}
$p_g(\calf_\xi)=h^{0,n}_B(\calf_\xi),$ the {\it arithmetic genus}
$p_a(\calf_\xi)=(-1)^n(\chi_{hol}(\calf_\xi)-1),$ and the irregularity
$q=q(\cals)=h^{0,1}_B.$ We remark that $q= {1\over 2}b^B_1(\calf_\xi)={1\over
2}b_1(M)$ is actually a topological invariant by (5) of
Proposition 1.9 of [BGN1]. In the case $n=2$, things
simplify much more. All the basic Hodge numbers are given in terms of the
three invariants $q,p_g(\calf_\xi)$ and $h^{1,1}_B(\calf_\xi),$ and we have
the relations
$$\chi_B(\calf_\xi)=2+2p_g(\calf_\xi)-4q+h^{1,1}_B(\calf_\xi), \qquad
\tau_B(\calf_\xi)=2+2p_g(\calf_\xi)-h^{1,1}_B(\calf_\xi) \leqno{\inv.7}$$ 
$$\chi_B(\calf_\xi)+\tau_B(\calf_\xi)=4\chi_{hol}(\calf_\xi).$$
The fact that $q$ is actually a topological invariant has some nice
consequences for $n=2.$

\noindent{\sc Proposition} \inv.8: \tensl Let $\cals=(\xi,\eta,\Phi,g)$ be a 
Sasakian structure on a 5-manifold $M$ with $b_1(M)=0.$ Then we have
$$b_2(\calf_\xi)=1+b_2(M), \qquad \chi_B(\calf_\xi)=3+b_2(M)\geq 3, \qquad
\chi_{hol}(\calf_\xi)=1+p_g(\calf_\xi)\geq 1.$$
\tenrm

\noindent{\sc Proof}: Since $b_1(M)=0$ Proposition 1.9 of [BGN1] and the exact
cohomology sequence \inv.5 imply
$$H^2_B(\calf_\xi)\approx \bbr\oplus H^2(M,\bbr), \leqno{\inv.9}$$
and the results easily follow from this and equations \inv.7. \hfill\za

\noindent{\sc Corollary} \inv.10: \tensl On a Sasakian 5-manifold $M$ with
$b_1(M)=0$ we must have $b_2(M)\geq 2p_g(\calf_\xi).$ In particular, every
Sasakian structure $\cals=(\xi,\eta,\Phi,g)$ on $S^5$ or $S^2\times S^3$
satisfies $p_g(\calf_\xi)=0$ and $\chi_{hol}(\calf_\xi)=1.$ \tenrm

\noindent{\sc Proof}: This follows from Proposition \inv.8 since by
Proposition 1.9 of [BGN1] we must have $h^{1,1}_B(\calf_\xi)\geq 1.$ \hfill\za

The reader no doubt notices that a classification of Sasakian structures on
compact 5-manifolds involves a ``transverse Enriques-Kodaira classification''
paralleling complex surface theory. Such a classification is currently under
study.

\bigskip
\baselineskip = 10 truept
\centerline{\bf \pos. Positive Sasakian Geometry}
\bigskip

\noindent{\sc Definition} \pos.1: \tensl A Sasakian manifold $M$ is said to
be {\it positive} if its basic first Chern class $c_1(\calf_\xi)$ can be
represented by a basic positive definite $(1,1)$-form. \tenrm

A basic ingredient in studying positive Sasakian geometry is the ``transverse
Yau Theorem'' of  El Kacimi-Alaoui [ElK]: 

\noindent{\sc Theorem} \pos.2 [ElK]: \tensl If $c_1(\calf_\xi)$ is represented
by a real basic $(1,1)$ form $\grr^T,$ then it is the Ricci curvature form
of a unique transverse K\"ahler form $\gro^T$ in the same basic cohomology
class as $d\eta.$  \tenrm

In [BGN1] this theorem was reformulated in terms of Sasakian geometry in a
convenient way, namely

\noindent{\sc Theorem} \pos.3 [BGN1]: \tensl Let $(M,\xi,\eta,\Phi,g)$ be a
Sasakian manifold whose basic first Chern class $c_1(\calf_\xi)$ is represented
by the real basic $(1,1)$ form $\grr,$ then there is a unique Sasakian
structure $(\xi,\eta_1,\Phi_1,g_1)\in\gF(\xi)$ homologous to
$(\xi,\eta,\Phi,g)$ such that $\grr_{g_1}=\grr-2d\eta_1$ is the Ricci form of
$g_1,$ and $\eta_1=\eta + \grz_1,$ with $\grz_1={1\over 2}d^c\phi.$ The metric
$g_1$ and endomorphism $\Phi_1$ are then given by 
$$\Phi_1=\Phi -\xi\otimes \grz_1\circ \Phi, \qquad
           g_1=d\eta_1\circ(\hbox{id}\otimes \Phi_1)+\eta_1\otimes \eta_1.$$
 \tenrm  

One key ingredient still missing in the study of transverse holomorphic
geometry is the transverse analogue of Dolbeault's Theorem. The reason is that
the sheaves in the basic Dolbeault complex are not fine, since they are
constant along the leaves of the foliation.  We can, however, obtain some
partial results.

\noindent{\sc Proposition} \pos.4: \tensl If $\cals=(\xi,\eta,\Phi,g)$ is a
complete (i.e. $g$ is complete) and positive Sasakian structure on a manifold
$M$ of dimension $2n+1,$ then $M$ is compact and $q=0.$ If in addition 
$\cals=(\xi,\eta,\Phi,g)$ is quasi-regular, then $h^{p,0}(\calf_\xi)=0$ for
all $p>0,$ and so $\chi_{hol}(\calf_\xi)=1, \quad
p_g(\calf_\xi)=p_a(\calf_\xi)=0.$  Moreover, any two quasi-regular positive
Sasakian structures on a 5-manifold must have the same basic Hodge number
$h^{1,1}_B$, the same basic Euler characteristic $\chi_B,$ and the same basic
Hirzebruch signature $\tau_B.$ \tenrm

\noindent{\sc Proof}: As mentioned above $q={1\over 2}b_1(M)$ follows from
Proposition 1.9 of [BGN1]. Since $c_1(\calf_\xi)>0$ Theorem
\pos.3 implies there is a Sasakian structure $(\xi,\eta_1,\Phi_1,g_1)\in\gF(\xi)$ 
homologous to $(\xi,\eta,\Phi,g)$ whose Ricci curvature is positive and
represents $c_1(\calf_\xi)$.  It then follows from Hasegawa and Seino's [HS]
application of Myers' Theorem to Sasakian geometry that  $g_1$ is complete, and
hence, $M$ is compact with finite fundamental group. This implies $q=0.$

To prove the second part we note that any nontrivial element of
$H^{p,0}_B(\calf_\xi)$ is represented by a basic $(p,0)$-form in the kernel of
$\bar{\partial},$ i.e. by a transversely holomorphic $p$-form $\gra.$ Since
$\cals=(\xi,\eta,\Phi,g)$ is quasi-regular, $M$ is the total space of an
orbifold $S^1$ V-bundle over a compact K\"ahler orbifold $\calz.$ Moreover,
since $\gra$ is basic, it descends to a nontrivial element $\check{\gra}\in
H^{p,0}(\calz).$ Now since the Sasakian structure $\cals$ is positive so is
the K\"ahler structure on $\calz,$ i.e. $c_1(\calz)>0,$ and this implies
$h^{p,0}(\calz)=0$ by the Kodaira-Baily vanishing Theorem. \hfill\za 

\noindent{\sc Proof of Theorem A}: $c_1(\calf_\xi)$ can be represented by a
positive definite basic $(1,1)$ form $\grr,$  so by Theorem \pos.3 there is a Sasakian
structure $\cals_1=(\xi,\eta_1,\Phi_1,g_1)\in\gF(\xi)$ homologous to
$(\xi,\eta,\Phi,g)$ such that $\grr_{g_1}=\grr-2d\eta_1$ is the Ricci form of
$g_1.$ Let $g^T_1$ denote the transverse K\"ahler metric of this Sasakian
structure. Then the Ricci curvatures of $g^T_1$ and $g_1$ are related by
$$\Ric_{g_1}|_{\cald_1\times \cald_1} =\Ric_{g^T_1} - 2g^T_1.
\leqno{\pos.3}$$ 
Next for any real number $a>0$ we can perform a homothetic
deformation [YK] of the Sasakian structure by defining
$$g^T_2={1\over a}g^T_1, \quad \eta_2={1\over a}\eta_1, \quad \xi_2=a\xi,
\quad \Phi_2=\Phi_1, \leqno{\pos.4}$$
in which case $\cals_2=(\xi_2,\eta_2,\Phi_2,g_2)$ is a Sasakian structure where
$g_2=g^T_2 +\eta_2\otimes \eta_2.$ Notice also that $\cals_1$ and $\cals_2$
both have the same contact subbundles with the same underlying transverse
complex structures, and the same characteristic foliations. Now since the
Ricci tensor is invariant under homothety we find 
$$\Ric_{g_2}|_{\cald_2\times \cald_2} =\Ric_{g^T_2} - {2\over a}g^T_2.
\leqno{\pos.5}$$
But since $Ric_{g^T_2}>0$ and $M$ is compact there exists $a_0\in \bbr^+$
such that for all $a>a_0$ we have $\Ric_{g_2}|_{\cald_2\times \cald_2}>0.$ But
also for any Sasakian metric we have $\Ric_{g_2}(X,\xi_2)=2n\eta_2(X)$ which
proves the result. \hfill\za

We shall make use of 

\noindent{\sc Proposition} \pos.6: \tensl Any simply connected compact positive
Sasakian manifold  $M$ is spin. \tenrm

\noindent{\sc Proof}: $M$ is spin if and only if its second Stiefel-Whitney
class $w_2(M)$ vanishes. But if $\cald$ is the contact subbundle of a Sasakian
structure $\cals=(\xi,\eta,\Phi,g)$ on $M,$ we have a natural splitting
$$TM=\cald \oplus L_\xi$$
where $L_\xi$ is the trivial real line bundle generated by $\xi.$ Thus, 
$$w_2(M)=w_2(TM)=w_2(\cald)$$
which is the mod 2 reduction of $c_1(\cald)\in H^2(M,\bbz).$ 

Now suppose that $\cals$ is a positive Sasakian
structure  so that $c_1(\calf_\xi)>0.$ Choose a transverse K\"ahler
metric $g'_T$ whose K\"ahler form $\gro'_T\in c_1(\calf_\xi).$ Then both
$\gro'_T$ and its transverse Ricci form $\grr'_T$ represent $c_1(\calf_\xi).$
This transverse K\"ahler structure defines a positive Sasakian structure
$\cals'=(\xi,\eta',\Phi',g')$ on $M$ such that $g'=g'_T+\eta'\otimes \eta',$
and $\gro'_T$ represents $c_1(\calf_\xi).$ 
The exact sequence 1.11 of [BGN1] becomes
$$\matrix{&H^2(M,\bbz)&\cr
          &\decdnar{}&\cr
          0\ra{1.2} \bbr\fract{\grd}{\ra{1.5}}
H^2_B(\calf_{\xi})\fract{\gri_*}{\ra{1.8}}&H^2(M,\bbr)&\ra{1.3} \cr}
\leqno{\pos.7}$$  
where $\grd(c)=c[\gro'_T]_B.$ But since $c_1(\calf_{\xi})$ is represented by
$\gro'_T$ we have $c_1(\cald')= \gri_*c_1(\calf_{\xi})=\gri_*\grd(1)=0$ in
$H^2(M,\bbr).$ But since $M$ is simply connected $c_1(\cald')=0,$ so this
implies that $w_2(M)=w_2(\cald')=0.$ \hfill\za

Actually we have proven more.

\noindent{\sc Proposition} \pos.8: \tensl For every positive Sasakian structure
$\cals=(\xi,\eta,\Phi,g)$ on $M$ the first Chern class $c_1(\cald)\in
H^2(M,\bbz)$ is torsion. \tenrm

So on a Sasakian manifold the free part of the first Chern class of the contact
bundle $\cald$ is an obstruction to positivity, and in particular to the
existence of a compatible Sasakian-Einstein metric. However, the authors know
of no example of a simply connected non-spin manifold admitting a Sasakian
structure, or more generally a Sasakian manifold with the free part of
$c_1(\cald)$ nontrivial. It should be mentioned that Sasakian manifolds are
known to admit $\hbox{Spin}^\bbc$-structures [Mor]. Now combining our results
with Smale's [Sm] well-known classification of simply connected spin
5-manifolds, we arrive at

\noindent{\sc Theorem} \pos.9: \tensl Let $M^5$ be a complete positive
Sasakian 5-manifold. Then the universal cover $\tilde{M^5}$ is diffeomorphic to
one of the following:
$$M^5(k,\gra_i)=S^5\# k\#(S^2\times S^3)\# M^5_{\gra_1}\#\cdots \#
M^5_{\gra_r}$$
for some nonnegative integers $k,r$ where $M^5_\gra$ is a compact
simply connected 5-manifold with $H_2(M^5_\gra,\bbz)\approx \bbz_\gra\oplus
\bbz_\gra,$ and $\gra_i$ are positive integers satisfying
$\gra_1|\gra_2|\cdots |\gra_r.$ \tenrm

The problem of existence of positive Sasakian structures on the above
5-manifolds is still open in general. However, in the next section we prove
existence in the case $r=0$ for all $k>0.$ We should mention that for $r=0$ and
$0\leq k \leq 9$ the manifolds are known to have \Se structures
[BG2, BGN1, BGN2] which are automatically positive. We should also mention
that generally there is an a priori obstruction for $M^{2n+1}$ to admit a
Sasakian structure, namely the top Stiefel-Whitney class $w_{2n+1}.$ However,
it follows from Smale's Theorem 3.2 [Sm] that $w_5(M^5(k,\gra_i))=0;$ hence,
all $M^5(k,\gra_i)$ are candidates for admitting Sasakian structures. These
5-manifolds with $\gra_i$ nontrivial can never be realized as hypersurfaces
in a well-formed weighted projective space since a result of [BG2] says that
such hypersurfaces necessarily have no torsion in $H_2(M,\bbz).$ Nevertheless,
other methods are available to show that $M^5(k,\gra_i)$ can admit Sasakian
structures. Such examples with $k=0$ and $r>0$ are given in [BL].

\bigskip
\baselineskip = 10 truept
\centerline{\bf \consum. Positive Sasakian Structures on $k\#(S^2\times S^3)$}
\bigskip

We consider weighted homogeneous polynomials $f$ of degree $d$ in four
complex variables $z_0,z_1,z_2,z_3$ which describe hypersurfaces in $\bbc^4$
with only an isolated singularity at the origin. Here we are interested in the
case that $f$ has  degree $d=k+1$ with weights $\bfw=(1,1,1,k),$ $k\geq 2$ and
Fano index $I=2.$ Explicitly $f$ is given by a general polynomial of the form
$$f(z_0,z_1,z_2,z_3)=g_{(k+1)}(z_0,z_1,z_2)+g_{(1)}(z_0,z_1,z_2)z_3
\leqno{\consum.1}$$ 
where $g_{(k+1)}(z_0,z_1,z_2)$ is a homogeneous polynomial of degree
$(k+1)$ which is not divisible by $g_{(1)}(z_0,z_1,z_2)\neq0$. 
This guarantees that $f$ has only an isolated
singularity at the origin. The zero locus of $f$ cuts out a cone $C_f\subset
\bbc^4$ and the link $L_f$ is defined by intersecting with the unit sphere
$S^7,$ viz. 
$$L_f=C_f\cap S^7. \leqno{\consum.2}$$
The link $L_f$ is a smooth simply connected 5-manifold which for $k\geq 5$
depends on $2(k-5)$ effective parameters. To see this, note that the group
$\gG(1,1,1,k)$ of complex automorphisms of the projective space
$\bbc\bbp(1,1,1,k)$ is obtained by the projectivisation of the following group
of transformations [BGN1]
$$\varphi_{\bba,\lambda,\phi}\pmatrix{z_0\cr z_1\cr z_2\cr z_3\cr}=
\pmatrix{\bba\pmatrix{z_0\cr z_1\cr z_2\cr}\cr
         \lambda z_3+\phi_{(k)}(z_0,z_1,z_2)\cr},$$
where $\bba\in GL(3,\bbc)$, $\lambda\in\bbc^*$, and 
$\phi_{(k)}(z_0,z_1,z_2)$ is an arbitrary homogeneous polynomial of degree 
$k$. 

To prove Theorem B of the introduction we give several lemmas.
\medskip
\noindent{\sc Lemma} \consum.2: \tensl The 5-manifold $L_f$ is spin. \tenrm

\noindent{\sc Proof}: We show that $L_f$ has natural positive Sasakian
structures. The result will then follow by Proposition \pos.6. It is shown in
[BG2] that links of isolated hypersurface singularities have natural Sasakian
structures $\cals$ whose characteristic foliation $\calf_\xi$ is determined by
the weights. Let $\calz_f$ denote the space of leaves of $\calf_\xi$ which is
embedded in the weighted projective space $\bbc\bbp(1,1,1,k).$ The Sasakian
structures will be positive if the basic first Chern class $c_1(\calf_\xi)$ is
positive. But since $c_1(\calf_\xi)$ is the pullback of the first Chern class
$c_1(\calz_f)$ of $\calz_f$ by the natural projection, it suffices to show
that $c_1(\calz_f)$ is positive. But 
$$c_1(\calz_f)\simeq c_1(K^{-1}_{\calz_f})\simeq c_1(\calo(2)_{\calz_f})$$ 
by 3.6 of [BGN1] which is clearly positive. \hfill\za
\medskip
\noindent{\sc Lemma} \consum.3: \tensl $b_2(L_f)=k.$ \tenrm

\noindent{\sc Proof}: The procedure for computing the Betti numbers of links
is given by Milnor and Orlik [MO]. Associate to
any monic polynomial $f$ with roots $\gra_1,\ldots,\gra_k\in \bbc^*$ its
divisor  
$$\hbox{div}~f= <\gra_1>+\cdots+<\gra_k>\leqno{\consum.4}$$
as an element of the integral ring $\bbz[\bbc^*]$ and let $\grL_n = \hbox{div}
(t^n-1).$  The `rational weights' used in [MO] are just ${d\over w_i},$
and are written in irreducible form, ${d\over w_i}={u_i\over v_i}$. The divisor
of the characteristic polynomial $\grD(t)$ is then determined by
$$\hbox{div}~\grD(t)~= \prod_i({\grL_{u_i}\over v_i}-1)~= 1+\sum
a_j\grL_j,\leqno{\consum.5.}$$
where $a_j\in \bbz$ and the second equality is obtained by using  the
relations  $\grL_a\grL_b=\break\gcd(a,b)\grL_{lcm(a,b)}.$ The second Betti
number of the link is then given by 
$$b_2(L_f)=1+\sum_ja_j.\leqno{\consum.6}$$

In our case we have
$$\hbox{div}~\grD(t)~=({\grL_{k+1}\over k}-1)(\grL_{k+1}-1)^3=
({\grL_{k+1}\over k}-1)(\grL_{k+1}-1)  ((k-1)\grL_{k+1}+1)=(k-1)\grL_{k+1}+1$$
implying $b_2(L_f)=k.$ \hfill\za

To prove Theorem B of the introduction it suffices to prove

\noindent{\sc Lemma} \consum.7: \tensl $L_f$ is diffeomorphic to
$k\#(S^2\times S^3).$ \tenrm

\noindent{\sc Proof}: This will follow from a Theorem of Smale [Sm] and Lemmas
\consum.2 and \consum.3 if we can show that $H_2(L_f,\bbz)$ has no torsion.
But since the gcd of any three weights is one this follows from Lemma 5.8 of
[BG2]. \hfill\za

Finally, we mention the failure of the sufficient conditions described in
[JK,BGN1] for $\calz_f$ to admit a K\"ahler-Einstein metric, and hence, for
$L_f$ to admit a \Se metric. Indeed, it follows from Lemma 5.1 of [BGN1] that
for any $k\geq 2,$ $(\calz_f,{2+\gre\over 3}D)$ cannot be Kawamata log terminal
for every effective $\bbq$-divisor that is numerically equivalent to
$K^{-1}_{\calz_f}.$ Furthermore, the Hitchin-Thorpe inequality ($c_1^2\geq 0$)
which prohibits  smooth compact complex surfaces with $q=p_g=0$ and $b_2\geq
9$ from admitting any Einstein metrics whatsoever, is ineffective for
hypersurfaces in the weighted projective space $\bbp(\bfw)$ since any such
hypersurface $\calz_f$ of degree $d$ satisfies $$c_1^2(\calz_f)=
{d(|\bfw|-d)^2\over w_0w_1w_2w_3}\geq 0.$$  
In our case we have $c_1^2=4(1+{1\over k}).$

\bigskip
\medskip
\centerline{\bf Bibliography}
\medskip
\font\ninesl=cmsl9
\font\bsc=cmcsc10 at 10truept
\parskip=1.5truept
\baselineskip=11truept
\ninerm

\item{[BG1]} {\bsc C. P. Boyer and  K. Galicki}, {\ninesl On Sasakian-Einstein
Geometry}, Int. J. of Math. 11 (2000), 873-909.
\item{[BG2]} {\bsc C. P. Boyer and  K. Galicki}, {\ninesl New Einstein Metrics
in Dimension Five},  math.DG/0003174, submitted for publication.
\item{[BGN1]} {\bsc C. P. Boyer, K. Galicki, and M. Nakamaye}, {\ninesl On the
Geometry of Sasakian-Einstein 5-Manifolds}, math.DG/0012047; submitted
for publication.
\item{[BGN2]} {\bsc C. P. Boyer, K. Galicki, and M. Nakamaye}, {\ninesl
Sasakian-Einstein Structures on $9\#(S^2\times S^3)$}, math.DG/0102181.
\item{[BL]} {\bsc V. Braun and C.-H. Lui}, {\ninesl On extremal transitions of
Calabi-Yau threefolds and the singularity of the associated 7-space from
rolling}, hep-th/9801175 v2.
\item{[DK]} {\bsc J.-P. Demailly and J. Koll\'ar}, {\ninesl Semi-continuity of
complex singularity exponents and K\"ahler-Einstein metrics on Fano
orbifolds}, preprint AG/9910118, to appear in Ann. Scient. Ec. Norm. Sup. Paris
\item{[ElK]} {\bsc A. El Kacimi-Alaoui}, {\ninesl Op\'erateurs transversalement
elliptiques sur un feuilletage riemannien et applications}, Compositio
Mathematica 79 (1990), 57-106.
\item{[ElKN]} {\bsc A. El Kacimi-Alaoui and M. Nicolau}, {\ninesl On the
topological invariance of the basic cohomology}, Math. Ann. 295 (1993),
627-634.
\item{[HS]} {\bsc I. Hasegawa and M. Seino}, {\ninesl Some remarks on Sasakian 
geometry--applications of Myers' theorem and the canonical affine connection}, 
J. Hokkaido Univ. Education 32 (1981), 1-7. 
\item{[JK]} {\bsc J.M. Johnson and J. Koll\'ar}, {\ninesl K\"ahler-Einstein
metrics on log del Pezzo surfaces in weighted projective 3-space}, preprint
AG/0008129, to appear in Ann. Inst. Fourier.
\item{[Mil]} {\bsc J. Milnor}, {\ninesl Singular Points of Complex
Hypersurfaces}, Ann. of Math. Stud. 61, Princeton Univ. Press, 1968.
\item{[MO]} {\bsc J. Milnor and P. Orlik}, {\ninesl Isolated singularities
defined by weighted homogeneous polynomials}, Topology 9 (1970), 385-393.
\item{[Mor]} {\bsc S. Moroianu}, {\ninesl Parallel and Killing spinors on
$\hbox{Spin}^c$-manifolds}, Commun. Math. Phys. 187 (1997), 417-427.
\item{[Sm]} {\bsc S. Smale}, {\ninesl On the structure of 5-manifolds},
Ann. Math. 75 (1962), 38-46.
\item{[SY]} {\bsc J.-P. Sha and D.-G Yang}, {\ninesl Positive Ricci
curvature on the connected sums of $S^n\times S^m$}, J. Diff. Geom. 33
(1991), 127-137.
\item{[Ton]} {\bsc Ph. Tondeur}, {\ninesl Geometry of Foliations}, Monographs
in Mathematics, Birkh\"auser, Boston, 1997.
\item{[YK]} {\bsc K. Yano and M. Kon}, {\ninesl
Structures on manifolds}, Series in Pure Mathematics 3,
World Scientific Pub. Co., Singapore, 1984.
\medskip
\bigskip \line{ Department of Mathematics and Statistics
\hfil March 2001} \line{ University of New Mexico \hfil }
\line{ Albuquerque, NM 87131 \hfil } \line{ email: cboyer@math.unm.edu,
galicki@math.unm.edu, nakamaye@math.unm.edu\hfil} \line{ web pages:
http://www.math.unm.edu/$\tilde{\phantom{o}}$cboyer,
http://www.math.unm.edu/$\tilde{\phantom{o}}$galicki \hfil}
\bye